\documentstyle[11pt]{article}
\begin{document}
\newcommand{\p}{\parallel }
\makeatletter \makeatother
\newtheorem{th}{Theorem}[section]
\newtheorem{lem}{Lemma}[section]
\newtheorem{de}{Definition}[section]
\newtheorem{rem}{Remark}[section]
\newtheorem{cor}{Corollary}[section]
\renewcommand{\theequation}{\thesection.\arabic {equation}}

\title{\bf The noncommutative family Atiyah-Patodi-Singer index theorem}

\author{ Yong Wang \\
}

\date{}
\maketitle

\begin{abstract} In this paper, we define the eta cochain form and prove its regularity when the kernel of a family of Dirac operators is a vector bundle.
We decompose the eta form as a pairing of the eta cochain form with the Chern character of an idempotent matrix and we also decompose the Chern character of the index bundle for a fibration with boundary as a pairing of the family Chern-Connes character for a manifold with boundary with the Chern character of an idempotent matrix. We define the family $b$-Chern-Connes character and then we prove that it is entire and
give its variation formula. By this variation formula, we prove another noncommutative family Atiyah-Patodi-Singer index theorem. Thus, we extend the results of Gezler and Wu to the family case.\\

\noindent{\bf Keywords:}\quad
 Eta cochain form; family Chern-Connes character for manifolds with boundary; family $b$-Chern-Connes character; variation formula.\\

\noindent{\bf MSC(2010):}\quad 58J20, 19K56\\
\end{abstract}

\section{Introduction}

\quad In [APS], Atiyah-Patodi-Singer introduced the eta invariant and proved their famous
Atiyah-Patodi-Singer index theorem for manifolds with boundary. In [BC], using Cheeger's cone method, Bismut and Cheeger defined the eta form which is a family version of the eta invariant and extended the APS index formula to the family case under the condition that all boundary Dirac operators are invertible. In [MP1,2], using the Melrose's $b$-calculus, Melrose
and Piazza extended the Bismut-Cheeger family index theorem to the case that boundary Dirac operators are not invertible.
In [Do], Donnelly extended the APS index theorem to the equivariant case by
modifying the Atiyah-Patodi-Singer original method. In [Zh], Zhang
got this equivariant Atiyah-Patodi-Singer index theorem by using a
direct geometric method in [LYZ].\\
\indent On the other hand, in [Wu], Wu proved the
Atiyah-Patodi-Singer index theorem in the framework of
noncommutative geometry. To do so, he introduced the eta
cochain (called the higher eta invariant in [Wu]) which is a
generalization of the classical Atiyah-Patodi-Singer eta
invariant in [APS], then proved its regularity by using the
Getzler symbol calculus [Ge1] as adopted in [BF] and
computed
 its radius of convergence. Subsequently, he proved the variation
formula of eta cochains, using which he got the noncommutative
Atiyah-Patodi-Singer index theorem. In [Ge2], using superconnection,
Getzler gave another proof of the noncommutative
Atiyah-Patodi-Singer index theorem, which was more difficult, but
avoided mention of the operators $b$ and $B$ of
cyclic cohomology. In [Wa1], we defined the equivariant eta cochain and proved its regularity using the method in [CH], [Fe] and [Zh]. Then
we proved an equivariant noncommutative
Atiyah-Patodi-Singer index theorem. In [Wa2], we defined infinitesimal equivariant eta cochains and proved their regularity. In [LMP], Lesch, Moscovici and Pflaum presented the Chern-Connes character of the Dirac operator associated to a
$b$-metric on a manifold with boundary in terms of a retracted cocycle in relative cyclic cohomology. Blowing-up the metric one recovered the pair of characteristic currents that
represent the corresponding de Rham relative homology class, while the blowdown
yielded a relative cocycle whose expression involves higher eta cochains and
their $b$-analogues. The corresponding pairing formula with relative $K$-theory
classes captured information about the boundary and allowed to derive geometric
consequences. In [Xi], Xie proved an analogue for odd dimensional manifolds with boundary,
in the $b$-calculus setting, of the higher Atiyah-Patodi-Singer index theorem
by Getzler and by Wu. Xie also obtained a natural counterpart of the eta invariant
for even dimensional closed manifolds.\\
\indent  The purpose of this paper is to extend the theorems due to Getzler and Wu to the family case. our main theorems are as follows (for related
definitions, see Sections 3-5).\\

\noindent {\bf Theorem 3.6}~{\it Suppose that all $D_{M,z}$ are invertible
with $\lambda$ the smallest positive eigenvalue of all $|D_{M,z}|$. We assume that
$||d(p|_M)||<\lambda$ and $p\in M_{r\times r}(C^{\infty}_*(N))$, then in the cohomology of $X$}
$${\rm ch}[{\rm Ind}(pD_{z,+,\varepsilon}p)]=\langle\tau(B),{\rm
Ch}(p)\rangle.\eqno(1.1)$$\\

\noindent {\bf Theorem 5.1} ~{\it  Suppose that all $D_{M,z}$ are invertible
with $\lambda$ the smallest positive eigenvalue of all $|D_{M,z}|$. We assume that
$||d(p|_M)||<\lambda$ and $p\in M_{r\times r}(C^{\infty}_{\rm exp}(\widehat{N}))$, then in the cohomology of $X$}
$${\rm ch}[{\rm Ind}(pD_{z,+}p)]=\int^b_{\widehat{N}/X}\widehat{A}(R^{{\widehat{N}}/X}){\rm ch}({\rm Imp})-\langle
{\widehat{\eta}}^*(B^M), {\rm ch}_*(p_M)\rangle.\eqno(1.2)$$

\indent The above theorems obviously apply with no mayor modifications to the twisted case by an extra odd differential form.
 This paper is organized as follows: in Section 2, we define
the eta cochain form and prove its regularity when the kernel of a family of Dirac operators is a vector bundle..
In Section 3,  we decompose the eta form as a pairing of the eta cochain form with the Chern character of an idempotent matrix and we also decompose the Chern character of the index bundle for a fibration with boundary as a pairing of the family Chern-Connes character for manifolds with boundary with the Chern character of an idempotent matrix. In Section 4, We define the family $b$-Chern-Connes character and then we prove that it is entire and
give its variation formula. In Section 5, by this variation formula, we prove another noncommutative family Atiyah-Patodi-Singer index theorem. Thus, we extend the results of Gezler and Wu to the family case.\\

\section{The eta cochain form }

\quad ~~In this Section, we define
the eta cochain form and prove its regularity.\\ \indent Firstly, we recall the Bismut superconnection.
 Let $M$ be a $n+{q}$ dimensional compact connected manifold without boundary
 and $X$ be a ${q}$ dimensional compact connected
 manifold without boundary. We assume that $\pi :M\rightarrow X$ is a submersion of
 $M$ onto $X$, which defines a fibration of $M$ with the fibre $Z$. For
 $y\in X$, $\pi^{-1}(y)$ is a submanifold $M_y$ of $M$. Denote by $TZ$
 the $n$-dimensional vector bundle on $M$ whose fibre $T_xM_{\pi x}$
 is the tangent space at $x$ to the fibre $M_{\pi (x)}$. We assume
 that $M$ and $X$ are oriented. We take a smooth
 horizontal subbundle $T^HM$ of $TM$. A vector field
 $X\in \Gamma (X,TX)$ will be identified with its horizontal lift $X^H\in  \Gamma (M, T^HM)$. Moreover $T^H_xM$ is isomorphic to $T_{\pi(x)}X$ via
 $\pi_*$. We take a Riemannian metric on $X$ and then lift the
 Euclidean scalar product $g_X$ of $TX$ to $T^HM$.
We further assume that $TZ$ is endowed with a scalar product $g_Z$. Thus
we can introduce on $TM$ a new scalar product $g_X\oplus g_Z$, and
denote by $\nabla^L$ the Levi-Civita connection on $TM$ with respect
to this metric. Set $\nabla^X$ denote the Levi-Civita connection on
$TX$ and we still denote by $\nabla^X$ the pullback connection on
$T^HM$. Let $\nabla^Z=P_Z(\nabla^L)$ where $P_Z$ denotes the orthogonal
projection to $TZ$. Set $\nabla^{\oplus}=\nabla^X\oplus \nabla^Z$
and $S=\nabla^L-\nabla^{\oplus}$ and $T$ be the torsion tensor of
$\nabla^{\oplus}$. Denote by $SO(TZ)$ the $SO(n)$ bundle of oriented
orthonormal frames in $TZ$. Now we assume that the bundle $TZ$ is spin.
Denote by $S(TZ)$ the associated spinor bundle and $\nabla^Z$ can be
lifted to a connection on $S(TZ)$. Let $D$ be the Dirac
operator in the tangent direction defined by $D=\sum_{j=1}^nc(e^*_j)\nabla^{S(TZ)}_{e_j}$ where $\nabla^{S(TZ)}$ is a spin connection on $S(TZ)$. Set $E$ be the vector bundle $\pi^*(\wedge T^*X)\otimes S(TZ)$. Then the Bismut superconnection acting on $E$ is defined by
$${{B}}=D+\sum_{\alpha=1}^{{q}}f_\alpha^*\wedge(\nabla^{S(TZ)}_{f_\alpha}+\frac{1}{2}k(f_\alpha))-\frac{1}{4}c(T),\eqno(2.1)$$
where
$$k(f_\alpha)=\sum_{j=1}^n\left<\nabla^{TZ}_{f_\alpha}e_j-[f_\alpha,e_j],e^j\right>,~~c(T)=-\sum_{\alpha<\beta}\sum_jf^\alpha\wedge f^\beta c(e_j)
\left<[f^H_\alpha,f^H_\beta],e_j\right>.\eqno(2.2)$$
\indent Let $\psi_t:dy_\alpha\rightarrow \frac{dy_\alpha}{\sqrt{t}}$ be the rescaling operator. Let ${{B}}_t=\sqrt{t}\psi_t({{B}})$
and $F_t={{B}}_t^2$. Let ${\rm tr^{even}}$ denote taking the trace with value in $\Omega^{\rm even}(X)$. When ${\rm dim}Z$ is odd, for $a_0,\cdots,a_{2k}\in C^{\infty}(M)$, we define
the family cochain ${\rm ch}_{2k}({{B}}_t,\frac{d{{B}}_t}{dt})$ by the formula:\\
$${\rm ch}_{2k}({{B}}_t,\frac{d{{B}}_t}{dt})(a_0,\cdots ,
a_{2k})~~~~~~~~~~~~~~~~~~~~~~$$
$$=
\sum^{2k}_{j=0}(-1)^j\langle a_0,[{{B}}_t,a_1],\cdots ,[{{B}}_t,a_{j}],\frac{d{{B}}_t}{dt},[{{B}}_t,a_{j+1}]
,\cdots , [{{B}}_t,a_{2k}]\rangle _t,\eqno(2.3)$$
 If $A_j~ (0\leq j\leq q)$ are operators on $\Gamma (E)$, we define:\\
\bigskip
$$\langle A_0,\cdots , A_q\rangle_t=\int_{\triangle_q}{\rm tr^{even}}[A_0e^{-\sigma_0F_t}A_1e^{-\sigma_1F_t}\cdots A_q
e^{-\sigma_qF_t}]d\sigma,\eqno(2.4)$$ where
$\triangle_q=\{(\sigma_0,\cdots,\sigma_q)|\sigma_0+\cdots+\sigma_q=1,~\sigma_j\geq 0\}$ is a simplex in ${\bf R^q}$.
When ${\rm dim}Z$ is even, in (2.4), we use ${\rm str}$ instead of ${\rm tr^{even}}$ and define ${\rm ch}_{2k}({{B}}_t,\frac{d{{B}}_t}{dt})$.\\
\indent We assume that the kernel of $D$ is a complex vector bundle. Formally, the {\bf eta cochain form} is defined to be an even cochain sequence
by the formula:\\
$$\widehat{\eta}_{2k}({{B}})=\frac{1}{\sqrt{\pi}}\int^{\infty}_{0}{\rm ch}_{2k}
({{B}}_t,\frac{d{{B}}_t}{dt})dt,~~{\rm when}~ {\rm dim}Z~~ {\rm is ~odd};
\eqno(2.5)$$
$$\widehat{\eta}_{2k}({{B}})=\int^{\infty}_{0}{\rm ch}_{2k}
({{B}}_t,\frac{d{{B}}_t}{dt})dt,~~{\rm when}~ {\rm dim}Z~~ {\rm is ~even}.
\eqno(2.6)$$
This integral makes sense by the following Lemma 2.1 and Lemma 2.5.
 Then $\widehat{\eta}_0({{B}})(1)$ is the eta form defined by Bismut and Cheeger in [BC].
In order to prove that the
above definition is well defined, it is necessary to check the
integrality near the two ends of the integration.
 Firstly, the regularity at infinity comes from the
following lemma.\\

 \noindent {\bf Lemma 2.1}~~{\it We assume that the kernel of $D$ is a complex vector bundle. For
$a_0,\cdots,a_{2k}\in C^{\infty}(M)$, we have}
$${\rm ch}_{2k}({{B}}_t,\frac{d{{B}}_t}{dt})(a_0,\cdots ,
a_{2k})=O(t^{-\frac{3}{2}}),~~{\rm as}~t\rightarrow\infty.\eqno(2.7)$$
\noindent{\bf {\it Proof.}}~~Since the kernel of $D$ is a complex vector bundle, our proof is very similar to the proof of Lemma 3.5 in [Wa2] (see revised version arXiv:1307.8189). We just use Lemma 9.4 in [BGV] instead of Lemma 3.4 in [Wa2]. We use $D+\frac{c(T)}{4}$ and $\psi_t:dy_\alpha\rightarrow \frac{dy_\alpha}{\sqrt{t}}$ instead of $D-\frac{c(\widehat{X})}{4}$ and $\psi_t:\widehat{X}\rightarrow \frac{\widehat{X}}{t}$ in Lemma 3.5 in [Wa2] respectively where $\widehat{X}$ is the Killing vector field. We note that ${\rm ch}_{2k}({{B}}_t,\frac{d{{B}}_t}{dt})$ corresponds to $\frac{1}{2\sqrt{t}}{\bf {\rm ch}}_k(\sqrt{t}D_{-\widehat{X}},D_{\widehat{X}})$ in [Wa2]. Comparing with the single operator case in Lemma 2 in [CM], the operator $[B_t,a_j]=\sqrt{t}[c(d_Za_j)+\frac{1}{\sqrt{t}}d_Xa_j\wedge]$ is instead
of $\sqrt{t}[D,a_j]$ and $\delta_t(g)$ in Lemma 9.21 in [BGV] emerges, where $\delta_t(g)=1+O(t^{-\frac{1}{2}})S_0$ and $S_0$ is a smooth operator. By these differences, in the discussions of Lemma 2 in [CM], the number of copies of $e^{-\sigma_ltD^2}(I-H)$ may be less than $\frac{k}{2}+1$. But the coefficients of $S_0$ and $d_Xa_j\wedge$ are $O(t^{-\frac{1}{2}})$. Through careful observations, we still get (2.7). ~~$\Box$\\

In the following, we prove the regularity at zero of the eta cochain form. We know that $\frac{d{{B}}_t}{dt}=\frac{1}{2\sqrt{t}}\psi_t
(D+\frac{c(T)}{4})$. We introduce the Grassmann variable $dt$ which anticommutates with $c(e_j)$ and $dy_\alpha$. Set $\widehat{F}=F+dt(D+\frac{c(T)}{4})$. Let
$${\rm ch}_{2k}(\widehat{F})(a_0,\cdots ,
a_{2k})
=t^k\int_{\triangle_{2k}}\psi_t{\rm tr^{even}}[a_0e^{-t\sigma_0\widehat{F}}[B,a_1]\cdots [B,a_{2k}]
e^{-t\sigma_{2k}\widehat{F}}]d\sigma.\eqno(2.8)$$
By the Duhamel principle and $(dt)^2=0$, we have
$$e^{-t\sigma_j\widehat{F}}=e^{-t\sigma_j{F}}-tdt\int^{\sigma_j}_0e^{-t(\sigma_j-\xi){F}}(D+\frac{c(T)}{4})e^{-t\xi{F}}d\xi.\eqno(2.9)$$
By (2.8) and (2.9), we have
$${\rm ch}_{2k}(\widehat{F})(a_0,\cdots ,
a_{2k})={\rm ch}_{2k}({F})(a_0,\cdots ,
a_{2k})+t^{\frac{3}{2}}{\rm ch}_{2k}({{B}}_t,\frac{d{{B}}_t}{dt})(a_0,\cdots ,
a_{2k})dt.\eqno(2.10)$$
Let ${A}$ be an
operator and $l$ be a positive interger.
 Write
$${A}^{[l]}=[\widehat{F},~{A}^{[l-1]}],~{A}^{[0]}=A,~ {A}^{(l)}=[F,{A}^{(l-1)}],~{A}^{(0)}=A.$$
Similar to Lemma 4.4 in [Wa3], we have\\

\noindent{\bf Lemma 2.2}~{\it Let $A$ a finite order fibrewise
differential
operator with form coefficients, then for any $s>0$, we have:}\\
$$e^{-sF}A=\sum^{N-1}_{l=0}\frac{(-1)^l}{l!}s^lA^{(l)}e^{-sF}+(-1)^Ns^NA^{(N)}(s);\eqno(2.11)$$
$$e^{-s\widehat{F}}A=\sum^{N-1}_{l=0}\frac{(-1)^l}{l!}s^lA^{[l]}e^{-s\widehat{F}}+(-1)^Ns^NA^{[N]}(s),\eqno(2.12)$$
\noindent {\it where $A^{(N)}(s)$ and $A^{[N]}(s)$ are given by}\\
$$A^{(N)}(s)=\int_{\triangle_N}e^{-u_1sF}A^{(N)}e^{-(1-u_1)sF}du_1du_2\cdots du_N;\eqno(2.13)$$
$$A^{[N]}(s)=\int_{\triangle_N}e^{-u_1s\widehat{{F}}}A^{[N]}e^{-(1-u_1)s\widehat{F}}du_1du_2\cdots du_N.\eqno(2.14)$$

\indent As in [CH], [Fe], [Wa3], by Lemma 2.2, we have for a sufficient large $N$,\\

\indent ${\rm ch}_{2k}(F)(a_0,\cdots,a_{2k})$
$$=\psi_t
\sum^{N}_{\lambda_1,\dots,\lambda_{2k}=0}
\frac{
(-1)^{\lambda_1+\cdots+\lambda_{2k}}}{\lambda_1!\cdots \lambda _{2k}!}Ct^{|\lambda|+{k}}
{\rm
tr^{even}}\left[a_0[B,a_1]^{(\lambda_1)}\cdots [B,a_{2k}]^{(\lambda_{2k})}e^{-tF}\right]+O(t^{\frac{3}{2}});\eqno(2.15)$$
\indent ${\rm ch}_{2k}(\widehat{F})(a_0,\cdots,a_{2k})$
$$=\psi_t
\sum^{N}_{\lambda_1,\dots,\lambda_{2k}=0}
\frac{
(-1)^{\lambda_1+\cdots+\lambda_{2k}}}{\lambda_1!\cdots \lambda _{2k}!}Ct^{|\lambda|+{k}}
{\rm
tr^{even}}\left[a_0[B,a_1]^{[\lambda_1]}\cdots [B,a_{2k}]^{[\lambda_{2k}]}e^{-t\widehat{F}}\right]+O(t^{\frac{3}{2}}),\eqno(2.16)$$
where $C$ is a constant. Recall Lemma 2.17 in [Wa3] which extends the corresponding Lemma in [Po] and [PW] \\
\indent Let $U$ be an open subset of ${\bf{R}}^n$. We define Volterra symbols and Volterra $\Psi DO$¡¯s on $U\times {\bf{R}}^{n+1}/0$
as follows.\\

\noindent {\bf {Definition 2.3}}~~{\it The set $S_V^m(U\times
{\bf{R}}^{n+1})\otimes \wedge T^*_zB,~m\in{\bf{Z}}$ ,
consists of smooth functions $q(x, \xi, \tau)$ on $U\times
{\bf{R}}^{n}\times {\bf{R}}$ with
an asymptotic expansion $q\sim \sum_{j\geq 0}q_{m-j},$ where:}\\
-{\it  $q_l\in C^{\infty}(U\times [({\bf{R}}^{n}\times {\bf{R}})\backslash 0])\otimes \wedge T^*_zB$ is a homogeneous Volterra symbol of degree $l$, i.e. $q_l$ is parabolic
homogeneous of degree $l$ and satisfies the property}

(i){\it ~ $q$ extends to a continuous function on $(
{\bf{R}}^n\times \overline{{\bf{C}}_-})\backslash 0 $ in such way to
be holomorphic in the
last variable when the latter is restricted to ${{\bf{C}}}_-$.}\
\\
- {\it The sign $\sim$ means that, for any integer $N$ and any compact $K\subset U,$ there is a constant
$C_{NK\alpha\beta k}>0$ such that for $x\in K$ and for $|\xi|+|\tau|^{\frac{1}{2}}>1$ we have}

$$||\partial_x^\alpha\partial_\xi^\beta\partial_\tau^k(q-\sum_{j<N}q_{m-j})(x,\xi,\tau)||\leq
C_{NK\alpha\beta k}(|\xi|+|\tau|^{\frac{1}{2}})^{m-N-|\beta|-2k}.\eqno(2.17)$$
{\it For $q=\sum_lq_l\omega^l$ where $q_l\in S_V^m(U\times
{\bf{R}}^{n+1})$ and $\omega^l\in \wedge^l T^*_zB$, we define
$||q||=\sum_l|q_l|||\omega^l||$ and $||\omega^l||$ is the norm of $\omega^l$ in  $(\wedge^l T^*_zB,g^{TB}_z).$}\\

\noindent {\bf Definition 2.4}  {\it The set $\Psi_V^m(U\times
{\bf{R}},\wedge T^*_zB), ~m\in{\bf{Z}}$ , consists of
continuous operators $Q$ from $C_c^{\infty}(U_x\times
{\bf{R}}_t,\wedge T^*_zB)$ to $C^{\infty}(U_x\times
{\bf{R}}_t,\wedge T^*_zB)$
such that:}\\
(i) {\it $ Q $ has the Volterra property;}\\
(ii) {\it $Q = q(x,D_x,D_t) + R$ for some symbol $q$ in $S_V^m(U\times {\bf{R}},\wedge T^*_zB)$
and some smoothing operator $R$.}\\

\indent In the sequel if $Q$ is a Volterra $\Psi DO$, we let $K_Q(x, y, t-s)$ denote its distribution kernel, so that
the distribution $K_Q(x, y, t)$ vanishes for $t< 0$.\\

\noindent{\bf Lemma 2.5}~ (Lemma 2.17 in [Wa3])~{\it Let $Q\in \Psi_V^*({\bf{R}}^n\times {\bf{R}}, S(T(M_z))\otimes \wedge ^*T^*_zX)$ have Getzler order $m$
and model operator $Q_{(m)}$. Then as $t\rightarrow 0^+$ we have:}
$$ 1)~\sigma[\psi_tK_Q(0,0,t)]^{(j)}=\omega^{{\rm odd}}O(t^{\frac{j-n-m-2}{2}})+O(t^{\frac{j-n-m-1}{2}}),~~{\rm if~} m-j~~{\rm ~is~ odd};$$
$2)~\sigma[\psi_tK_Q(0,0,t)]^{(j)}=t^{\frac{j-n-m-2}{2}}K_{Q_{(m)}}(0,0,1)^{(j)}+\omega^{{\rm odd}}O(t^{\frac{j-n-m-1}{2}})+O(t^{\frac{j-n-m}{2}}),$
$~~~~~~~~~~{\rm if~} m-j~~{\rm ~is~ even},$\\
{\it where $[K_Q(0,0,t)]^{(j)}$ denotes the degree $j$ form component in $M_z$ and $\omega^{{\rm odd}}O(t^{\frac{j-n-m-2}{2}})$ denotes
that the coefficients of $t^{\frac{j-n-m-2}{2}}$ are in $\wedge^{{\rm odd}}(T^*X)\otimes \wedge(T^*(M_z))$.}\\

\noindent {\bf Lemma 2.6} {\it The following estimate holds}\\
$${\rm ch}_{2k}({{B}}_t,\frac{d{{B}}_t}{dt})\sim O(1) ~~{\rm when}~t\rightarrow 0,\eqno(2.18)$$\\

\noindent {\bf Proof.} By (2.10), (2.15) and (2.16), in order to prove Lemma 2.6, we only prove
$$\psi_t
t^{|\lambda|+{k}}
{\rm
tr^{even}}[a_0[B,a_1]^{[\lambda_1]}\cdots [B,a_{2k}]^{[\lambda_{2k}]}e^{-t\widehat{F}}]$$
$$
-\psi_t
t^{|\lambda|+{k}}
{\rm
tr^{even}}[a_0[B,a_1]^{(\lambda_1)}\cdots [B,a_{2k}]^{(\lambda_{2k})}e^{-tF}]=O(t^{\frac{3}{2}})dt.\eqno(2.19)$$
This a local problem and we fix a point $x_0$ in $M_z$. Set
 $$h(x)=1+\frac{1}{2}dt\sum_{j=1}^nx_jc(e_j)\eqno(2.20)$$ as in [Zh].
By (5.29) in [Wa3], we have
$$h[F+dt(D+\frac{c(T)}{4})]h^{-1}=F+dtu,\eqno(2.21)$$
where the Getzler order $O_G(u)\leq 0$ of $u$. Write
$$\widetilde{A}^{[l]}=[h\widehat{F}h^{-1},~\widetilde{A}^{[l-1]}],~\widetilde{A}^{[0]}=A.$$
Then
$$\psi_t
t^{|\lambda|+{k}}
{\rm
tr^{even}}[a_0[B,a_1]^{[\lambda_1]}\cdots [B,a_{2k}]^{[\lambda_{2k}]}e^{-t\widehat{F}}]$$
$$
=\psi_t
t^{|\lambda|+{k}}
{\rm
tr^{even}}[a_0\widetilde{[B,a_1]}^{[\lambda_1]}\cdots \widetilde{[B,a_{2k}]}^{[\lambda_{2k}]}e^{-th\widehat{F}h^{-1}}].\eqno(2.22)$$
By the Volterra calculus, we have
$$(\frac{\partial}{\partial t}+F+dtu)^{-1}=(\frac{\partial}{\partial t}+F)^{-1}-dt(\frac{\partial}{\partial t}+F)^{-1}u(\frac{\partial}{\partial t}+F)^{-1}.\eqno(2.23)$$
Let
$$a_0\widetilde{[B,a_1]}^{[\lambda_1]}\cdots \widetilde{[B,a_{2k}]}^{[\lambda_{2k}]}
=A_0+dtA_1,\eqno(2.24)$$
where
$$A_0=a_0{[B,a_1]}^{(\lambda_1)}\cdots {[B,a_{2k}]}^{(\lambda_{2k})}.$$
Then
$$a_0\widetilde{[B,a_1]}^{[\lambda_1]}\cdots \widetilde{[B,a_{2k}]}^{[\lambda_{2k}]}(\frac{\partial}{\partial t}+F+dtu)^{-1}
-A_0(\frac{\partial}{\partial t}+F)^{-1}$$
$$=-A_0dt(\frac{\partial}{\partial t}+F)^{-1}u(\frac{\partial}{\partial t}+F)^{-1}+dtA_1(\frac{\partial}{\partial t}+F)^{-1}.\eqno(2.25)$$
By (2.24) in [Wa3], in order to prove (2.19), we only need to prove
$$t^{k+|\lambda|}\psi_t{\rm tr^{even}}[A_0(\frac{\partial}{\partial t}+F)^{-1}u(\frac{\partial}{\partial t}+F)^{-1}]=O(t^{\frac{3}{2}}),\eqno(2.26)$$
$$t^{k+|\lambda|}\psi_t{\rm tr^{even}}[A_1(\frac{\partial}{\partial t}+F)^{-1}]=O(t^{\frac{3}{2}}).\eqno(2.27)$$
We note that $A_0(\frac{\partial}{\partial t}+F)^{-1}u(\frac{\partial}{\partial t}+F)^{-1}\in {\rm End}^-(\wedge^*(TX)\otimes S(TZ))$, so
when we take ${\rm tr^{even}}$, only the coefficient of $c(e_1)\cdots c(e_n)$ is left and other terms are zero. Note that
$$O_G( t^{k+|\lambda|}A_0(\frac{\partial}{\partial t}+F)^{-1}u(\frac{\partial}{\partial t}+F)^{-1})\leq -4,\eqno(2.28)$$
so by Lemma 2.5 (1) for $j=n$ odd and $m=-4$ and taking ${\rm tr^{even}}$, we get (2.26). By $O_G(u)\leq 0$ and (2.24), we get
$O_G(t^{|\lambda|+k}A_1(\frac{\partial}{\partial t}+F)^{-1})\leq -4$. Again $j=n$, so we get (2.27). Thus we prove Lemma 2.6.~~$\Box$\\

\noindent{\bf Remark.} We also introduce a new Bismut superconnection on $\widetilde{M}=M\times {\bf{R}}_+ \rightarrow X\times {\bf{R}}_+$ as in [BGV, Thm. 10.32] and prove a formula which is similar to (2.10). Then we can give a new proof of Lemma 2.6 as in [BGV, p. 347].\\

\indent For the idempotent $p\in {\cal M}_r(C^{\infty}(M))$, its
Chern character ${\rm Ch}(p)$ in entire cyclic homology is defined
by the formula (for more details see [GS]):
$${\rm Ch}(p)={\rm Tr}(p)+\sum_{k\geq 1}\frac{(-1)^k(2k)!}{k!}{\rm
Tr}_{2k}((p-\frac{1}{2})\otimes \overline{p}^{\otimes 2k})
\eqno(2.29)$$ where
$${\rm Tr}_{2k}:~{\cal M}_r(C^{\infty}(M))\otimes\left({\cal
M}_r(C^{\infty}(M))/{\cal M}_r({\bf C})\right)^{\otimes
2k}\rightarrow C^{\infty}(M)\otimes(C^{\infty}(M)/{\bf C})^{\otimes
2k}$$ is the generalized trace map. Let
$$||dp||=||[B,p]||=\sum_{i,j}||d_Mp_{i,j}||\eqno(2.30)$$
where $p_{i,j}~(1\leq i,j\leq r)$ is the entry of $p$. Similar to Proposition 2.17 in [Wa1], we have\\

\noindent {\bf Proposition 2.7} {\it Suppose that all $D_z$ are invertible
with $\lambda$ the smallest positive eigenvalue of all $|D_z|$. We assume that
$||dp||<\lambda$, then
the pairing $\langle\widehat{\eta}^*({{B}}),{\rm Ch}_*(p)\rangle$ is well-defined.}\\

\section{The family index pairing for manifolds with boundary}

\quad In this section, we decompose the eta form as a pairing of the eta cochain form with the Chern character of an idempotent matrix and we also decompose the Chern character of the index bundle for a fibration with boundary as a pairing of the family Chern-Connes character for a manifold with boundary with the Chern character of an idempotent matrix.\\
\indent  Suppose that all $D_z$ are invertible
with $\lambda$ the smallest positive eigenvalue of all $|D_z|$ and
$||dp||<\lambda$.  Let $H=\Gamma(M,\wedge^*(TX)\otimes S(TZ))$ and
$$p(B\otimes I_r)p:~p(H\otimes {\bf C^r})=L^2(M,\wedge^*(TX)\otimes S(TZ)\otimes
p({\bf C^r}))$$
$$\rightarrow L^2(M,\wedge^*(TX)\otimes S(TZ)\otimes p({\bf C^r}))$$ be the Bismut superconnection with the coefficient from $F=p({\bf C^r})$.
Then we have\\

\noindent{\bf Theorem 3.1} {\it Under the assumption as above, we have up to an exact form on $X$
 $$\widehat{\eta}(p(B\otimes
 I_r)p)=\langle\widehat{\eta}^*({{B}}),{\rm Ch}_*(p)\rangle,\eqno(3.1)$$
where the left term is the Bismut-Cheeger eta form.}\\

 \indent Let
 $${\bf B}=\left[\begin{array}{lcr}
  \ 0 &  -B\otimes I_r \\
    \  B\otimes I_r  & 0
\end{array}\right];~
{\bf p}=\left[\begin{array}{lcr}
  \ p &  0 \\
    \ 0  & p
\end{array}\right];~
\sigma=\sqrt{-1}\left[\begin{array}{lcr}
  \ 0 &  I_r \\
    \  I_r  & 0
\end{array}\right]
$$
be operators from $H\otimes {\bf C^r}\oplus H\otimes {\bf C^r}$ to
itself, then
$${\bf B}\sigma=-\sigma {\bf B};~~\sigma {\bf p}={\bf p}\sigma.\eqno(3.2)$$
Moreover ${\bf B}e^{t{\bf B}^2}$ and $e^{t{\bf B}^2}~ (t>0)$ are
traceclass. For $u\in [0,1]$, let
$$B_u=(1-u)B+u[pBp+(1-p)B(1-p)]=B+u(2p-1)[B,p],\eqno(3.3)$$
then
$${\bf
B}_u=\left[\begin{array}{lcr}
  \  0 &  -B_u \\
    \ B_u  & 0
\end{array}\right]
={\bf B}+u(2p-1)[{\bf B},p].\eqno(3.4)$$ \noindent We consider the infinite dimensional bundle $H\otimes {\bf C^r}\oplus H\otimes {\bf C^r}$ on
 $X\times [0,1]\times {\bf R}\times [0,\infty)$, parameterized by
$(b,u,s,t)$. Let
$$\widetilde{{\bf B}}=t^{\frac{1}{2}}\psi_t{\bf
B}_u+s\sigma({\bf p}-\frac{1}{2}),\eqno(3.5)$$ then $A=d_{(u,s,t)}+\widetilde{{\bf B}}$ be a
superconnection on $H\otimes {\bf C^r}\oplus H\otimes {\bf C^r}$. Direct computations show that
$$(d+\widetilde{{\bf B}})^2=t\psi_t{\bf
B}_u^2-s^2/4-(1-u)t^{\frac{1}{2}}s\sigma[\psi_t{\bf B},p]+ds\sigma
(p-\frac{1}{2})$$
$$+t^{\frac{1}{2}}du(2p-1)[\psi_t{\bf
B},p]+\frac{1}{2}t^{-\frac{1}{2}}dt\psi_t[{\bf D}_u+\frac{c(T)}{4}].\eqno(3.6)$$
We also consider $A$ as $A_t$, which is a family superconnection
parameterized  by $t$ on the superbundle with the base $X\times [0,1]\times
{\bf R}$ and the fibre $H\otimes {\bf C^r}\oplus H\otimes {\bf C^r}$.
\noindent Let $\Gamma_u=\{u\}\times {\bf R}\subset [0,1]\times {\bf
 R}$ be a contour oriented in the direction of increasing $s$
 and $\gamma_s=[0,1]\times \{s\}$
 be a contour oriented in the direction of increasing $u$ . By the
Duhamel principle and the
 Stokes theorem as in page 225 in [Wa1], then
 $$d_X\omega=\int_{[0,1]\times {\bf
 R}}d\int_0^{+\infty}{\rm
 Str^{even}}(e^{A^2})=\left(\int_{\Gamma_1}-\int_{\Gamma_0}
-\int_{\gamma_{+\infty}}+\int_{\gamma_{-\infty}}\right)
\left[\int_0^{+\infty}{\rm
 Str^{even}}(e^{A^2})\right],\eqno(3.7)$$
 where ${\rm
 Str^{even}}$ denotes taking the supertrace with value in $\Omega^{\rm even}(X)\otimes \Omega([0,1]\times
{\bf R})$.
So in the cohomology of $X$, we have
 $$\int_{\Gamma_0}\int_0^{+\infty}{\rm
 Str^{even}}(e^{A^2})=\int_{\Gamma_1}\int_0^{+\infty}{\rm
 Str^{even}}(e^{A^2}).\eqno(3.8)$$
Similar to (3.8) in [Wa1], we have
$$\int_{\Gamma_0}\int_0^{+\infty}{\rm
 Str^{even}}(e^{A^2})=-4\sqrt{-1}\pi[\langle \widehat{\eta}^*({{B}}),{\rm
Ch}(p)\rangle-\langle \widehat{\eta}^*({{B}}),{\rm rk}(p){\rm
Ch}_*(1)\rangle].\eqno(3.9)$$
Similar to (3.10) in [Wa1], we have
$$\int_{\Gamma_1}\int_0^{+\infty}{\rm
 Str^{even}}(e^{A^2})=-2\sqrt{-1}\int_{-\infty}^{+\infty}e^{-s^2/4}ds$$
 $$\cdot
\int_0^{+\infty}\psi_t{\rm
 Tr^{even}}[(p-\frac{1}{2})
(D_1+\frac{c(T)}{4})e^{-t B_1^2}]d\sqrt{t}.\eqno(3.10)$$
By the following lemma 3.2 and (3.8)-(3.10), similar to (3.12) in [Wa1], we can prove Theorem 3.1.~~$\Box$\\

\noindent{\bf Lemma 3.2}~{\it Let $B_s=B+s(2p-1)[B,p]$ for $s\in [0,1]$. We assume that all $D_z$ be invertible and $||d_Mp||<\lambda$, then we have}
$\widehat{\eta}(B_0)=\widehat{\eta}(B_1)$.\\

\noindent{\bf Proof.} By $||d_Mp||<\lambda$, then $D_s=D+s(2p-1)[D,p]$ is invertible for $s\in [0,1]$. Similar to the discussions of Proposition 4.4 in [Wu], the eta form of $B_s$ is well defined. So $\widehat{\eta}(B_s)$ is smooth. Let $B_s=D_s+A_{[1]}
-\frac{c(T)}{4}$ and $A_0=(2p-1)[D,p]$. Then by the definition of the eta form and the Duhamel principle, we have
$$\frac{d}{ds}\widehat{\eta}(B_s)=\frac{1}{\sqrt{\pi}}\int^{+\infty}_0\psi_t{\rm tr^{even}}[A_0e^{-tB^2_s}]d\sqrt{t}+L,\eqno(3.11)$$
where $$
L=-\frac{1}{\sqrt{\pi}}\int^{+\infty}_0\psi_t{\rm tr^{even}}\left\{t(B_s+\frac{c(T)}{2}-A_{[1]})\right.$$
$$\left.\cdot\int^1_0e^{-\sigma tB^2_s}[(2p-1)[B,p],B_s]_+e^{-(1-\sigma) tB^2_s}d\sigma \right\}d\sqrt{t}.\eqno(3.12)$$
By ${\rm tr^{even}}(AB)={\rm tr^{even}}(BA)$ and $B_se^{-\sigma B^2_s}=e^{-\sigma B^2_s}B_s$, we have
$${\rm tr^{even}}\left\{B_s\int^1_0e^{-\sigma tB^2_s}[(2p-1)[B,p],B_s]_+e^{-(1-\sigma) tB^2_s}d\sigma\right\}$$
$$
=\int^1_0{\rm tr^{even}}\left\{(2p-1)[B,p]e^{-\sigma tB^2_s}[B_s,B_s]_+e^{-(1-\sigma) tB^2_s}d\sigma\right\},\eqno(3.13)$$
and
$${\rm tr^{even}}\left\{(\frac{c(T)}{2}-A_{[1]})\int^1_0e^{-\sigma tB^2_s}[(2p-1)[B,p],B_s]_+e^{-(1-\sigma) tB^2_s}d\sigma\right\}$$
$$=\int^1_0{\rm tr^{even}}\left\{(2p-1)[B,p]e^{-\sigma tB^2_s}[\frac{c(T)}{2}-A_{[1]},B_s]_+e^{-(1-\sigma) tB^2_s}d\sigma\right\}.\eqno(3.14)$$
By (3.12)-(3.14)
$$L=-\int^{+\infty}_0\frac{\sqrt{t}}{2\sqrt{\pi}}\psi_t\int^1_0{\rm tr^{even}}\left\{(2p-1)[B,p]e^{-\sigma tB^2_s}[D_s+\frac{c(T)}{4},B_s]_+e^{-(1-\sigma) tB^2_s}d\sigma\right\}dt,\eqno(3.15)$$
By $$\frac{d(t\psi_tB^2_s)}{dt}=\frac{1}{2}\psi_t[D_s+\frac{c(T)}{4},B_s]_+,\eqno(3.16)$$
(3.11) and (3.15), using the Duhamel principle and the Leibniz rule, then we get
$$\frac{\partial}{\partial t}\left\{\frac{\sqrt{t}}{\sqrt{\pi }}\psi_t{\rm tr^{even}}\left[(2p-1)[D,p]e^{-tB^2_s}\right]\right\}=
\frac{\partial}{\partial s}\left\{\frac{1}{2\sqrt{\pi t}}\psi_t{\rm tr^{even}}\left[(D_s+\frac{c(T)}{4})e^{-tB^2_s}\right]\right\}.\eqno(3.17)$$
So
$$\frac{d}{ds}\widehat{\eta}(B_s)=\frac{\sqrt{t}}{\sqrt{\pi }}\psi_t{\rm tr^{even}}\left[(2p-1)[D,p]e^{-tB^2_s}\right]|^{+\infty}_0.\eqno(3.18)$$
By $D_s$ being invertible, ${\rm tr^{even}}\left[(2p-1)[D,p]e^{-tB^2_s}\right]$ exponentially decays, so
$${\rm lim}_{t\rightarrow +\infty}\frac{\sqrt{t}}{\sqrt{\pi }}\psi_t{\rm tr^{even}}\left[(2p-1)[D,p]e^{-tB^2_s}\right]=0.\eqno(3.19)$$
By Lemma 2.3, similar to the discussions on page 164 in [Wu], we have
$${\rm lim}_{t\rightarrow 0}\frac{\sqrt{t}}{\sqrt{\pi }}\psi_t{\rm tr^{even}}\left[(2p-1)[D,p]e^{-tB^2_s}\right]$$
$$=c_0\int_Z\widehat{A}(TZ){\rm tr}\left\{(2p-1)(d_Zp){\rm exp}[\frac{\sqrt{-1}}{2\pi}(A'\wedge A'+dA')]\right\}=0,\eqno(3.20)$$
 where $A'=s(2p-1)d_Mp$. Then by (3.18)-(3.20), we prove Lemma 3.2.~~$\Box$\\

\indent Let $N$ be a fibration with the even-dimensional compact spin fibre. Let $M$ be the boundary of $N$. We endow $N$ with a metric which is a product in
a collar neighborhood of $M$. Denote by $B~(B_M)$ the Bismut superconnection on $N~(M)$. Let $C^{\infty}_*(N)=\{f\in C^{\infty}(N)|f$ is
independent of the normal coordinate $x_n$ near the boundary
$\}.$\\

 \noindent {\bf Definition
3.3}~The family Chern-Connes character on $N$, $\tau=\{
\tau_0,\tau_2,\cdots, \tau_{2q}\cdots \}$ is defined by
$$ \tau_{2q}(B)(f^0,f^1,\cdot,f^{2q}):= -\widehat{\eta}_{2q}(B_M)
(f^0|_M,f^1|_M,\cdot,f^{2q}|_M)$$
$$+\frac{1}{(2q)!(2\pi\sqrt{-1})^q}
\int_Z \widehat{A}(TZ)f^0 df^1\wedge \cdots\wedge
df^{2q},\eqno(3.21)$$
where $f^0,f^1,\cdot,f^{2q}\in C^{\infty}_*(N)$.\\

\indent Similar to Proposition 4.2 in [Wa1], we have\\

 \noindent {\bf Proposition 3.4}~{\it The family Chern-Connes character is
$b-\widetilde{B}$ closed (Here we use $\widetilde{B}$ instead of the Connes operator $B$. For the definitions of $b,~\widetilde{B}$, see [Co]). That is, in the cohomology of $X$, we have}
$$b\tau_{2q-2}+\widetilde{B}\tau_{2q}=0.\eqno(3.22) $$\\

\indent By Proposition 2,7, we have\\

\noindent {\bf Proposition 3.5} {\it Suppose that all $D_{M,z}$ are invertible
with $\lambda$ the smallest positive eigenvalue of all $|D_{M,z}|$. We assume that
$||d(p|_M)||<\lambda$, then
the pairing $\langle\tau,{\rm Ch}(p)\rangle$ is well-defined.}\\

\indent We let $C_1(M)=M\times (0,1],~\widetilde{N}=N\cup_{M\times \{1\}}C_1(M)$
and $\cal{U}$ be a collar neighborhood of $M$ in $N$. For
$\varepsilon>0$, we take a metric $g^{\varepsilon}$ of $\widetilde{N}$ such that
on ${\cal{U}}\cup_{M\times \{1\}}C_1(M)$
$$g^{\varepsilon}=\frac{dr^2}{\varepsilon}+r^2g^{M}.$$
\noindent Let $S=S^+\oplus S^-$ be spinors bundle associated to
$(\widetilde{N}_z,g^{\varepsilon})$ and $H^{\infty}$ be the set
$\{\xi\in\Gamma(\widetilde{N}_z,S)|~\xi ~{\rm and~ its~ derivatives~ are~ zero~
near~ the~ vertex~ of~ cone~}\}.$ Denote by $L^2_c(\widetilde{N}_z,S)$ the
$L^2-$completion of $H^{\infty}$ (similar define $L^2_{c}(\widetilde{N}_z,S^+)$
and $L^2_{c}(\widetilde{N}_z,S^-)$). Let
$$D_{z,\varepsilon}:~H^{\infty}\rightarrow H^{\infty};
~~D_{z,+,\varepsilon}:~H_+^{\infty}\rightarrow H_-^{\infty},$$ be the
Dirac operators associated to $(\widetilde{N}_z,g^{\varepsilon})$ which are
Fredholm operators for the sufficient small $\varepsilon.$ When $D_{M_z}$ is invertible, the index bundle of $\{D_z\}$ is well defined by [BC].
We recall the Bismut-Cheeger family index theorem for the twisting bundle ${\rm Im}p$ with the connection $pd$ in [BC]
$${\rm ch}[{\rm Ind}(pD_{z,+,\varepsilon}p)]=\sum_{r=0}^{\infty} \frac{(-1)^r}{r!(2\pi\sqrt{-1})^r}
\int_{Z}
\widehat{A}(TZ){\rm
Tr}[p(dp)^{2r}]-\widehat{\eta}(pB_Mp).\eqno(3.22)$$

Then we get\\

\noindent {\bf Theorem 3.6}~{\it Suppose that all $D_{M,z}$ are invertible
with $\lambda$ the smallest positive eigenvalue of all $|D_{M,z}|$. We assume that
$||d(p|_M)||<\lambda$ and $p\in M_{r\times r}(C^{\infty}_*(N))$, then in the cohomology of $X$}
$${\rm ch}[{\rm Ind}(pD_{z,+,\varepsilon}p)]=\langle\tau(B),{\rm
Ch}(p)\rangle.\eqno(3.23)$$\\

\noindent {\it Proof.}
 Let $\widehat{\tau}_{2q} $ be defined by
$$\widehat{\tau}_{2q}(B)(f^0,f^1,\cdot,f^{2q}):=
\frac{1}{(2q)!(2\pi\sqrt{-1})^q}
\int_Z \widehat{A}(TZ)f^0 df^1\wedge \cdots\wedge
df^{2q},\eqno(3.24)$$
where $f^0,f^1,\cdot,f^{2q}\in C^{\infty}_*(N)$. Recall
$${\rm ch}({\rm Im}p)=\sum_{q=0}^{\infty}\frac{(-1)^q}{(2\pi\sqrt{-1}q!}{\rm Tr}[p(dp)^{2q}],\eqno(3.25)$$
By (3.22) and Theorem 3.1 and (2.29), (3.24) and ${\rm Tr}[(dp)^{2k}]=0$ for $1\leq k$,
we have
$$\left<\widehat{\tau}_{\star}(B),{\rm ch}(p)\right>=\int_Z\widehat{A}(TZ){\rm ch}({\rm Im}p),\eqno(3.26)$$
so Theorem 3.6 holds.\\  $\Box$

\section{The family $b$-Chern-Connes character}

\quad In this section, we define the family $b$-Chern-Connes character which is the family version of the Getzler's $b$-Chern-Connes character
in [Ge2] and then we prove that it is entire and
give its variation formula.\\
\indent Let us recall the exact $b$-geometry (see [LMP],[Xi]). Let $N$ be a compact fibration with boundary $M$ and denote by $N^\circ$ its interior of $N$. We take the $b$-metric
$g_b=\frac{1}{r^2}dr\otimes dr+g_M$ near the $M$ where $r$ is the normal coordinate near the boundary. Let $x=lnr$ which gives an isometry
between the infinite cylinder $((-\infty,c_0]\times M,g_{cyl}=dx\otimes dx+g_M)$ and the collar neighborhood $U$ with the exact $b$-metric. Now we consider the complete Riemannian manifold $\widehat{N}=(-\infty,c_0]\times M\cup_M(N\backslash{U^\circ})$ instead of $N^\circ$ with the exact $b$-metric. Let
$C^{\infty}_{\rm exp}(\widehat{N})$ be the space of smooth functions on $\widehat{N}$ which expands exponentially on
the infinite cylinder $(-\infty,c_0]\times M$. A smooth function $f \in C^{\infty}(\widehat{N})$
expands exponentially on $(-\infty,c_0]\times M$ if
$f(x,y)\sim\sum_{k=0}^{\infty}e^{kx}f_k(y)$ for $(x,y)\in (-\infty,c_0]\times M$, where $f_k(y)\in C^{\infty}(M)$ for each $k$. That is
$$f(x,y)-\sum_{k=0}^{N-1}e^{kx}f_k(y)=e^{Nx}R_N(x,y),\eqno(4.1)$$
where all derivative of $R_N(x,y)$ in $x$ and $y$ are bounded. \\
\indent On $(-\infty,c_0]\times M$, we write $a=a_c+e^xa_{\infty}$ for $a\in C^{\infty}_{\rm exp}(\widehat{N})$ with $a_c,a_{\infty}\in C^{\infty}(\widehat{N})$ and $a_c$ constant with respect to $x$. Following [Xi],
define the $b$-norm of $a$ by $^b||a||:=||a_c||_1+2||a_{\infty}||_1$. The $b$-integral of $a$ along the fibre is defined by
$$\int^bad{\rm vol}:=\int_{N_z\backslash{U_z^\circ}}a|_{N_z\backslash{U_z^\circ}}d{\rm vol}+\int_{(-\infty,c_0]\times M_z}e^xa_{\infty}d{\rm vol}.\eqno(4.2)$$
Following [LMP, A.1] and [MP1], we can define the $b$-pseudodifferential operator with coefficients in $\wedge^*(TX)$ and
the pointwise trace of the Schwartz kernel of smooth $b$-pseudodifferential operators is a $b$-function. We define the $b$-trace is the $b$-integral of
this $b$-function. That is, for $A\in \Psi^{-\infty}_b(\widehat{N_z},\wedge^*(TX)\otimes S(T\widehat{N}_z))$ and its Schwartz kernel $k_A$,
define the $b$-trace which is in $\Omega(X)$ by
$$^b{\rm Str}(A)=\int^b{\rm Str}(k_A(x,x))d{\rm vol}.\eqno(4.3)$$
Let $B$ be the Bismut superconnection on $\widehat{N}$ and $B_t=t\psi_tB$ and $F_t=B^2_t$. By [MP1], $e^{-F_t}\in
\Psi^{-\infty}_b(\widehat{N_z},\wedge^*(TX)\otimes S(T\widehat{N}_z))$. For $A_0,\cdots,A_q\in \Psi^{\infty}_b(\widehat{N},\wedge^*(TX)\times S(T\widehat{N}_z))$,
we define
$$\langle\langle A_0,\cdots , A_q\rangle\rangle_b=\int_{\triangle_q}{\rm ^bStr}[A_0e^{-\sigma_0F}A_1e^{-\sigma_1F}\cdots A_q
e^{-\sigma_qF}]d\sigma,\eqno(4.4)$$
and
$$\langle\langle A_0,\cdots , A_q\rangle\rangle_{b,t}=\int_{\triangle_q}{\rm ^bStr}[A_0e^{-\sigma_0F_t}A_1e^{-\sigma_1F_t}\cdots A_q
e^{-\sigma_qF_t}]d\sigma.\eqno(4.5)$$
For $f_0,\cdots,f_k\in C^{\infty}_{\rm exp}(\widehat{N})$, we define the {\bf family $b$-Chern-Connes character} by
$$^b{\rm ch}^k(B)(f_0,\cdots,f_k):=\langle\langle f_0,[B,f_1],\cdots ,[B,f_k]\rangle\rangle_b;\eqno(4.6)$$
$$^b{\rm ch}^k(B_t)(f_0,\cdots,f_k):=\langle\langle f_0,[B_t,f_1],\cdots ,[B_t,f_k]\rangle\rangle_{b,t}.\eqno(4.7)$$
Define
$$^b{\rm ch}^k(B,V):=\sum_{0\leq j\leq k}(-1)^{j{\rm deg}V}\langle\langle f_0,[B,f_1],\cdots,[B,f_j],V,[B,f_{j+1}],\cdots,[B,f_k]\rangle\rangle_{b}.\eqno(4.8)$$
Similarly we may define $^b{\rm ch}^k(B_t,V)$. The family $b$-Chern-Connes character is well defined by the following Proposition 4.7. We recall the following lemma\\

\noindent{\bf Lemma 4.1} ([MP1, Proposition 9]) {\it For $A\in \Psi^{\infty}_{b,cl}(\widehat{N},\wedge^*(TX)\times S(T\widehat{N}_z))$
and  $L\in \Psi^{-\infty}_{b,cl}(\widehat{N},\wedge^*(TX)\times S(T\widehat{N}_z))$, we have}
$$^b{\rm Str}[A,L]=\frac{\sqrt{-1}}{2\pi}\int^{+\infty}_{-\infty}{\rm Str}_{M}\left(\frac{\partial I(A,\lambda)}{\partial \lambda}\cdot I(L,\lambda)\right)d\lambda,\eqno(4.9)$$
{\it where $I(L,\lambda)$ is the indicial family of $L$ (for the definition, see [LMP] or [MP1]).}\\

\indent Let ${\bf D}$ be the Dirac operator on the cylinder $(-\infty,+\infty)\times M$, then ${\bf D}=c(dx)\frac{d}{dx}+{\bf D_\partial}$.
On the boundary, $c(dx)$ gives a natrual identification of the even and odd half spinor bundle, then with respect to this splitting
$${\bf D}=\left(\begin{array}{lcr}
  \ 0 & -1 \\
    \ 1  & 0
\end{array}\right)
\frac{d}{dx}+
\left(\begin{array}{lcr}
  \ 0 & D_\partial \\
    \ D_\partial  & 0
\end{array}\right)
.\eqno(4.10)$$
By [MP1, p.139], we have\\

\noindent{\bf Lemma 4.2} {\it The following equality holds}
$$I(B,\lambda)=\sqrt{-1}c(dx)\lambda+B'^{M};~~I(F,\lambda)=\lambda^2+(B'^{M})^2\eqno(4.11)$$
{\it where with respect the above splitting}
$$B'^{M}={D_\partial}\left(\begin{array}{lcr}
  \ 0 & 1 \\
    \ 1  & 0
\end{array}\right)+\sum_{\alpha=1}^{{q}}f_\alpha^*\wedge(\nabla^{S(TM/X)}_{f_\alpha} +\frac{1}{2}k^M(f_\alpha))\left(\begin{array}{lcr}
  \ 1 & 0 \\
    \ 0  & 1
\end{array}\right)$$
$$~~~~~~~~~~~~~~~~~~~~~~~~~-\frac{1}{4}c(T^M)\left(\begin{array}{lcr}
  \ 0 & 1 \\
    \ 1  & 0
\end{array}\right).\eqno(4.12)$$\\

\indent By Lemma 4.2, we have
$$F'^M:=(B'^{M})^2\in \Omega^{\rm even}(X)\left(\begin{array}{lcr}
  \ L_1 & 0 \\
    \ 0  & L_1
\end{array}\right)+
\Omega^{\rm odd}(X)\left(\begin{array}{lcr}
  \ 0 & L_2 \\
    \ L_2  & 0
\end{array}\right),\eqno(4.13)$$
where $L_1,L_2 \in {\rm End}(S(TM_z))$. Similarly, we have
$$I([B,a],\lambda)=\left(\begin{array}{lcr}
  \ 0 & [D_\partial,a_\partial] \\
    \ [D_\partial,a_\partial]  & 0
\end{array}\right)+
d_Xa_\partial
\left(\begin{array}{lcr}
  \ 1 & 0 \\
    \ 0  & 1
\end{array}\right).\eqno(4.14)$$
\indent By Lemma 4.1 and Lemma 4.2, we have\\

\noindent{\bf Lemma 4.3}~ {\it For $K\in \Psi^{-\infty}_{b,cl}(\widehat{N},\wedge^*(TX)\times S(T\widehat{N}_z))$, we have}
$$^b{\rm Str}[B,K]=d_X {\rm ^bStr}(K)-\frac{1}{2\pi}\int^{+\infty}_{-\infty}{\rm Str}_M[c(dx)I(K,\lambda)]d\lambda.\eqno(4.15)$$

By Lemmas 4.1-4.3, similar to Lemma 6.3 in [Ge2], we have\\

\noindent{\bf Lemma 4.4}~{\it Let $A_j\in \Psi^{\infty}_{b,cl}(\widehat{N},\wedge^*(TX)\times S(T\widehat{N}_z))$ which indicial family is independent of
$\lambda$.\\
1. If $\varepsilon_j=(|A_0|+\cdots+|A_{j-1}|)(|A_j|+\cdots+|A_{k}|)$, then
$$\langle\langle A_0,\cdots,A_k\rangle\rangle_{b,t}=(-1)^{\varepsilon_j}\langle\langle A_j,\cdots,A_k,A_0,\cdots,A_{j-1}\rangle\rangle_{b,t}.\eqno(4.16)$$
2.
$$\langle\langle A_0,\cdots,A_k\rangle\rangle_{b,t}=\sum_{j=0}^k(-1)^{\varepsilon_j}\langle\langle 1, A_j,\cdots,A_k,A_0,\cdots,A_{j-1}\rangle\rangle_{b,t}.\eqno(4.17)$$
3.
$$-d_X\langle\langle A_0,\cdots,A_k\rangle\rangle_{b,t}+\sum_{j=0}^k(-1)^{|A_0|+\cdots+|A_{j-1}|}\langle\langle A_0,\cdots,[B_t,A_j],\cdots,A_k \rangle\rangle_{b,t}$$
$$=\langle\langle (A_0)_\partial,\cdots,(A_k)_\partial\rangle\rangle_{\partial,t},\eqno(4.18)$$
where when ${\rm dim}M_z$ is odd,
$$\langle\langle (A_0)_\partial,\cdots,(A_k)_\partial\rangle\rangle_{\partial,t}:=
\frac{-1}{2\sqrt{\pi}}\int_{\triangle_k}{\rm Str}_M\left[c(dx)A_{0,\partial}e^{-\sigma_0{F'_t}^M}\cdots A_{k,\partial}e^{-\sigma_k{F'_t}^M}\right]d\sigma.\eqno(4.19)$$
4. For $0\leq j<k$,
$$\langle\langle A_0,\cdots,[F_t,A_j],\cdots,A_k\rangle\rangle_{b,t}~~~~~~~~~~~~~~~~~~~~~~~~~~~~~~~~~~~~~~~~~~$$
$$=-\langle\langle A_0,\cdots,A_jA_{j+1},\cdots,A_k\rangle\rangle_{b,t}+
\langle\langle A_0,\cdots,A_{j-1}A_j,\cdots,A_k\rangle\rangle_{b,t}.\eqno(4.20)$$
For $j=k$,}
$$\langle\langle A_0,\cdots,A_{k-1},[F_t,A_k]\rangle\rangle_{b,t}~~~~~~~~~~~~~~~~~~~~~~~~~~~~~~~~~~~~~~~~$$
$$=\langle\langle A_0,\cdots,A_{k-2},A_{k-1}A_{k}\rangle\rangle_{b,t}-
(-1)^{\varepsilon_k}\langle\langle A_k A_0,A_1,\cdots,A_{k-1}\rangle\rangle_{b,t}.\eqno(4.21)$$\\

\noindent {\bf Proof.} 1)  By the definition of trace, we have
$$\langle\langle A_0,\cdots,A_k\rangle\rangle_{b,t}-(-1)^{\varepsilon_j}\langle\langle A_j,\cdots,A_k,A_0,\cdots,A_{j-1}\rangle\rangle_{b,t}$$
$$=\int_{\triangle_k} {\rm ^bstr}\left[A_0e^{-\sigma_0F_t}A_1\cdots A_{j-1}e^{-\sigma_{j-1}F_t}, A_je^{-\sigma_{j}F_t}A_{j+1}\cdots A_k e^{-\sigma_kF_t}\right]d\sigma.\eqno(4.22)$$
By $I(F_t,\lambda)=t\lambda^2+F^\partial_t$ and Lemma 4.1 and $\int_{-\infty}^{+\infty}\lambda e^{-t\lambda^2}d\lambda=0$, we know that
1) holds.\\
(2) (2) comes from (1) by the same trick in the [Get2, p.18].\\
(3)By $B_te^{F_t}=e^{F_t}B_t$, we have
$$\sum_{j=0}^k(-1)^{|A_0|+\cdots+|A_{j-1}|}\langle\langle A_0,\cdots,[B_t,A_j],\cdots,A_k \rangle\rangle_{b,t}$$
$$=^b{\rm Str}[B_t,A_0e^{-s_1F_t}A_1e^{-(s_2-s_1)F_t}\cdots A_ke^{-(1-s_k)F_t}].\eqno(4.23)$$
Then by Lemma 4.3 and (4.23),we get (3).\\
(4) By the Duhamel principle, we have
$$\frac{d}{ds_j}[e^{-(s_j-s_{j-1})F_t}A_je^{-(s_{j+1}-s_{j})F_t}]=-e^{-(s_j-s_{j-1})F_t}[F_t,A_j]e^{-(s_{j+1}-s_{j})F_t}.\eqno(4.24)$$
By the integration along $\int^{s_{j+1}}_{s_{j-1}}$, we have
$$[e^{-(s_{j+1}-s_{j-1})F_t},A_j]=-\int^{s_{j+1}}_{s_{j-1}}e^{-(s_j-s_{j-1})F_t}[F_t,A_j]e^{-(s_{j+1}-s_{j})F_t}ds_j.\eqno(4.25)$$
By (4.25), we get (4.20). By (4.25) and the indicial family of $A_k$ is independent of
$\lambda$ and Lemma 4.1, we get (4.21).  $\Box$

\indent By Lemma 4.4, similar to the proof of Theorem 6.2 in [Ge2], we get\\

\noindent{\bf Theorem 4.5} {\it When ${\rm dim}M_z$ is odd, for any $k\geq 0$, the following equality holds}
$$b{\rm ^bch}^{k-2}(B_t)+\widetilde{B}{\rm ^bch}^k(B_t)-d_X{\rm ^bch}^{k-1}(B_t)=\widetilde{Ch}^{k-1}({B'_t}^M)\circ {i_M}^*,\eqno(4.26)$$
{\it where }
$$\widetilde{Ch}^{k-1}({B'_t}^M)\circ i^*_M(a_0,\cdots,a_{k-1})=\langle\langle (a_0)_\partial,[{B'_t}^M,a_{1,\partial}],\cdots,[{B'_t}^M,a_{k-1,\partial}]\rangle\rangle_{\partial,t}.\eqno(4.27)$$\\

\noindent{\bf Proof.} By Lemma 4.4 (3), for $A_0=a_0,~A_j=[B_t,a_j]$ and $1\leq j\leq k-1$ we have
$$-d_X\langle\langle a_0,[B_t,a_1],\cdots,[B_t,a_{k-1}]\rangle\rangle_{b,t}+
+\langle\langle [B_t,a_0],[B_t,a_1],\cdots,[B_t,a_{k-1}]\rangle\rangle_{b,t}$$
$$+
\sum_{j=1}^{k-1}(-1)^{j-1}\langle\langle a_0,[B_t,a_1],\cdots,[B^2_t,a_j],\cdots,[B_t,a_{k-1}]\rangle\rangle_{b,t}$$
$$=\langle\langle (a_0)_\partial,[{B'_t}^M,a_{1,\partial}],\cdots,[{B'_t}^M,a_{k-1,\partial}]\rangle\rangle_{\partial,t}.\eqno(4.28)$$
By the definition of $\widetilde{B}$ (see [Co]) and Lemma 4.4 (2), we get
$$\widetilde{B}{\rm ^bch}^k(B_t)(a_0,\cdots,a_{k-1})=\langle\langle [B_t,a_0],[B_t,a_1],\cdots,[B_t,a_{k-1}]\rangle\rangle_{b,t}.\eqno(4.29)$$
By Lemma 4.4 (4), and $[B_t, a_ja_{j+1}]=[B_t, a_j]a_{j+1}+a_j[B_t,a_{j+1}]$, we have
$$b{\rm ^bch}^{k-2}(B_t)(a_0,\cdots,a_{k-1})=\sum_{j=1}^{k-1}(-1)^{j-1}\langle\langle a_0,[B_t,a_1],\cdots,[B^2_t,a_j],\cdots,[B_t,a_{k-1}]\rangle\rangle_{b,t}.\eqno(4.30)$$
By (4.28)-(4.30), we get Theorem 4.5. $\Box$\\

\indent By Theorem 4.5, we have\\

\noindent{\bf Theorem 4.6}~{\it When ${\rm dim}N_z$ is even and $k-1$ is even, the following equality holds}
$$\frac{d{\rm ^bch}^{k-1}(B_t)}{dt}+b{\rm ^bch}^{k-2}(B_t,\frac{dB_t}{dt})
+\widetilde{B}{\rm ^bch}^{k}(B_t,\frac{dB_t}{dt})$$
$$+d_X{\rm ^bch}^{k-1}(B_t,\frac{dB_t}{dt})
=-\frac{1}{\sqrt{\pi}}{\rm ch}^{k-1}(B^M_t,\frac{dB^M_t}{dt})
.\eqno(4.31)$$\\

\noindent{\bf Proof.} We know that $B_t$ is a superconnection on the infinite dimensional bundle $C^{\infty}(N,E)\rightarrow X$ which
we write ${\mathcal{E}}\rightarrow X$. Let $\widetilde{X}=X\times {\bf{R}}_+$, and $\widetilde{\mathcal{E}}$ be the superbundle
$\pi^*{\mathcal{E}}$ over $\widetilde{X}$, which is the pull-back to $\widetilde{X}$ of ${\mathcal{E}}$. Define a superconnection
$\widehat{B}$ on $\widetilde{\mathcal{E}}$ by the formula
$$(\widehat{B}\beta)(y,t)=(B_t\beta(\cdot,t))(y)+dt\wedge\frac{\partial\beta(y,t)}{\partial t}.\eqno(4.32)$$
The curvature $\widehat{\mathcal{F}}$ of $\widehat{B}$ is
$$\widehat{\mathcal{F}}={\mathcal{F}}_t-\frac{dB_t}{dt}\wedge dt,\eqno(4.33)$$
where ${\mathcal{F}}_t=B^2_t$ is the curvature of $B_t$. By the Duhamel principle, then
$$e^{-\widehat{\mathcal{F}}}=e^{-{\mathcal{F}}_t}-dt\left(\int^1_0e^{-u{\mathcal{F}}_t}\frac{dB_t}{dt}e^{-(1-u){\mathcal{F}}_t}du\right).\eqno(4.34)$$
Then for any $l\geq 0$, we have
$${\rm ^bch}^l(\widehat{B})={\rm  ^bch}^l(B_t)-dt{\rm ^bch}^l(B_t,\frac{dB_t}{dt}).\eqno(4.35)$$
By Theorem 4.5, we have
$$b{\rm ^bch}^{k-2}(\widehat{B})+\widetilde{B}{\rm ^bch}^k(\widehat{B})-d_X{\rm ^bch}^{k-1}(\widehat{B})=\widetilde{Ch}^{k-1}({\widehat{B'}}^M).\eqno(4.36)$$
By Theorem 4.5 and (4.36), (4.35) and $d_{\widetilde{X}}=d_X+dt\frac{d}{dt}$, we have
$$dt\left[\frac{d{\rm ^bch}^{k-1}(B_t)}{dt}+b{\rm ^bch}^{k-2}(B_t,\frac{dB_t}{dt})
+\widetilde{B}{\rm ^bch}^{k}(B_t,\frac{dB_t}{dt})\right.$$
$$\left.+d_X{\rm ^bch}^{k-1}(B_t,\frac{dB_t}{dt})\right]=\widetilde{Ch}^{k-1}({B'_t}^M)-\widetilde{Ch}^{k-1}({\widehat{B'}}^M).\eqno(4.37)$$
By (4.19), (4.35), (4.12) and (4.13), we get
$$\widetilde{Ch}^{k-1}({B'_t}^M)-\widetilde{Ch}^{k-1}({\widehat{B'}^M})=-\frac{1}{\sqrt{\pi}}dt{\rm ch}^{k-1}(B^M_t,\frac{dB^M_t}{dt})
.\eqno(4.38)$$
By (4.37) and (4.38), we get (4.31).~~$\Box$\\

\indent We recall that an even cochain $\{\Phi_{2n}\}$ is called
entire if $\sum_n||\Phi_{2n}||n!z^n$ is entire, where $||\Phi||:={\rm sup}_{^b||f^j||\leq 1}\{|\Phi(f^0,f^1,\cdots,f^{2k})|\}$ for $f_j\in C^{\infty}_{\rm exp}(\widehat{N})$.
 Then we have\\

\noindent{\bf Proposition 4.7} ~{\it ${\rm ^bch}(B)$ is an entire cochain and $\langle {\rm ^bch}(B),{\rm ch}(p)\rangle$ is well defined.}\\

\noindent{\bf Proof.}
For $A\in \Psi^{-\infty}_b(\widehat{N_z},\wedge^*(TX)\otimes S(T\widehat{N}_z))$ and its Schwartz kernel $k_A$, we define
$${\rm Str}^{N\setminus U}(A)=\int_{N_z\setminus U_z}{\rm Str}(k_A(x,x))d{\rm vol};~~^b{\rm Str^{end}}(A)=\int^b_{(-\infty,c_0)\times M_z}{\rm Str}(k_A(x,x))d{\rm vol}.\eqno(4.39)$$
So for $a_0,\cdots,a_q\in C^{\infty}_{\rm exp}(\widehat{N})$,
$$\int_{\triangle_q}{\rm ^bStr}\left[a_0e^{-\sigma_0F}[B,a_1]e^{-\sigma_1F}\cdots [B,a_q]
e^{-\sigma_qF}\right]d\sigma$$
$$=\int_{\triangle_q}{\rm Str}^{N\setminus U}\left[a_0e^{-\sigma_0F}[B,a_1]e^{-\sigma_1F}\cdots [B,a_q]
e^{-\sigma_qF}\right]d\sigma$$
$$+
\int_{\triangle_q}{\rm ^bStr^{end}}\left[a_0e^{-\sigma_0F}[B,a_1]e^{-\sigma_1F}\cdots [B,a_q]
e^{-\sigma_qF}\right]d\sigma.\eqno(4.40)$$
By the discussions on the compact fibration as in [BeC], we have
$$\left|\int_{\triangle_q}{\rm Str}^{N\setminus U}\left[a_0e^{-\sigma_0F}[B,a_1]e^{-\sigma_1F}\cdots [B,a_q]
e^{-\sigma_qF}\right]d\sigma\right|\leq {\rm Tr}(e^{-\frac{D^2}{2}})^b||a_0||^b||a_1||\cdots^b||a_q||.\eqno(4.41)$$
On the cylinder, we get
$$[B,a_j]=C_j+e^xB_j;~~a_0=C_0+e^xB_0,\eqno(4.42)$$
 where $$C_j=c(d_{M_z}(a_j)_c)+d_X(a_j)_c;~~B_j=c(d_{N_z}(a_j)_{\infty})+c((a_j)_{\infty}dx)+d_X(a_j)_\infty,$$
and $C_j$ is constant along the normal direction $x$. The second term in (4.40) can be written as a sum of terms of the following two types:\\
~ I)~$\int_{\triangle_q}{\rm ^bStr^{end}}\left[C_0e^{-\sigma_0F}C_1e^{-\sigma_1F}\cdots C_q
e^{-\sigma_qF}\right]d\sigma,$\\
~ II)~$\int_{\triangle_q}{\rm ^bStr^{end}}\left[A_0e^{-\sigma_0F}\cdots e^{-\sigma_jF}e^xB_je^{-\sigma_{j+1}F}\cdots A_q
e^{-\sigma_qF}\right]d\sigma,$ where $A_j=C_j~{\rm or} ~e^xB_j.$\\
Firstly we estimate the type I) integral. Without generality, we set $q=1$. Let $B^2=D^2+A_{[+]}$ and $D_R$ be the Dirac operator on the cylinder
$(-\infty,c_0)\times M_z$. By the Duhamel principle, we have
$$C_0e^{-\sigma_0F}C_1e^{-\sigma_1F}=C_0\sum_{m\geq 0}(-\sigma_0)^m\int_{\triangle_m}e^{-\sigma_0v_0D^2}A_{[+]}\cdots A_{[+]}e^{-\sigma_0v_mD^2}dv$$
$$\times C_1\sum_{l\geq 0}(-\sigma_1)^l\int_{\triangle_l}e^{-\sigma_1v'_0D^2}A_{[+]}\cdots A_{[+]}e^{-\sigma_1v'_lD^2}dv'$$
$$=C_0\sum_{m\geq 0}(-\sigma_0)^m\int_{\triangle_m}[e^{-\sigma_0v_0D^2}-e^{-\sigma_0v_0D_R^2}]
A_{[+]}\cdots A_{[+]}e^{-\sigma_0v_mD^2}dv$$
$$
\times C_1\sum_{l\geq 0}(-\sigma_1)^l\int_{\triangle_l}e^{-\sigma_1v'_0D^2}A_{[+]}\cdots A_{[+]}e^{-\sigma_1v'_lD^2}dv'$$
$$+\cdots+
C_0\sum_{m\geq 0}(-\sigma_0)^m\int_{\triangle_m}e^{-\sigma_0v_0D_R^2}A_{[+]}\cdots A_{[+]}e^{-\sigma_0v_mD^2_R}dv$$
$$\times C_1\sum_{l\geq 0}(-\sigma_1)^l\int_{\triangle_l}e^{-\sigma_1v'_0D^2_R}A_{[+]}\cdots A_{[+]}[e^{-\sigma_1v'_lD^2}-e^{-\sigma_1v'_lD^2_R}]dv'$$
$$+C_0\sum_{m\geq 0}(-\sigma_0)^m\int_{\triangle_m}
e^{-\sigma_0v_0D_R^2}A_{[+]}\cdots A_{[+]}e^{-\sigma_0v_mD^2_R}dv$$
$$\times C_1\sum_{l\geq 0}(-\sigma_1)^l\int_{\triangle_l}e^{-\sigma_1v'_0D^2_R}A_{[+]}\cdots A_{[+]}e^{-\sigma_1v'_lD^2_R}dv'.\eqno(4.43)$$
We know that $A_{[+]}$ is independent of $x$ on the cylinder and $D^2_R=\triangle_R+D^2_{M_z}$, so
$${\rm ^bStr^{end}}\left[C_0\sum_{m\geq 0}(-\sigma_0)^m\int_{\triangle_m}
e^{-\sigma_0v_0D_R^2}A_{[+]}\cdots A_{[+]}e^{-\sigma_0v_mD^2_R}dv\right.$$
$$\times C_1\sum_{l\geq 0}(-\sigma_1)^l\int_{\triangle_l}e^{-\sigma_1v'_0D^2_R}A_{[+]}\cdots A_{[+]}e^{-\sigma_1v'_lD^2_R}dv'\left|_{(-\infty,c_0)\times M_z}\right]=0.\eqno(4.44)$$
We estimate the first term $K_1$ in (4.43) and the estimate of other terms is similar. Since $D$ and $D_R$ are self adjoint, we can apply the functional calculus to these two operators. Then $||e^{-uD^2}||\leq 1$ and $||e^{-uD^2_R}||\leq 1$ for $u\geq 0$. By Theorem 3.2 (1) in [LMP],
similar to the proof of Lemma 2.2 in [Wa2], then $||K_1||_1$ is bounded. By the measure of the boundary of the simplex being zero,
we can estimate $K_1$ in the interior of the simplex, that is $\sigma_0>0,\sigma_1>0,v_j>0,v'_j>0$. We note that the zero order $b$-pseudodifferential operator is bounded and
$$||(1+D^2)^{-\frac{1}{2}}e^{-uD^2}||\leq L_0u^{-\frac{1}{2}};~~||e^{-uD^2}-e^{-uD_R^2}||_1\leq L'_0 u^r,\eqno(4.45)$$
where $L_0,L'_0$ are constant and $r$ is any integer. It holds that (see line 7 in [BC, P.21])
$$\int_{\triangle_m}v_0^{-\frac{1}{2}}\cdots v_{m-1}^{-\frac{1}{2}}dv=\frac{\pi^{\frac{m}{2}}}{\frac{m}{2}\Gamma(\frac{m+1}{2})}.\eqno(4.46)$$
When $m$ is odd, then
$$\frac{\pi^{\frac{m}{2}}}{\frac{m}{2}\Gamma(\frac{m+1}{2})}\leq \frac{2\pi^{\frac{m}{2}}}{(\frac{m+1}{2})!}.\eqno(4.47)$$
When $m$ is even, then
$$\frac{\pi^{\frac{m}{2}}}{\frac{m}{2}\Gamma(\frac{m+1}{2})}\leq \frac{2\pi^{\frac{m}{2}}}{(\frac{m}{2})!}.\eqno(4.48)$$
 So by (4.45)-(4.48), we get
\begin{eqnarray*}
||K_1||&\leq& ||C_0||||C_1||\sum_{l.m}\int_{\triangle_1}(\sigma_0)^m(\sigma_1)^l\int_{\triangle_m}
||e^{-\sigma_0v_0D^2}-e^{-\sigma_0v_0D_R^2}||_1\\
&&
\cdot||A_{[+]}(1+D^2)^{-\frac{1}{2}}||||(1+D^2)^{\frac{1}{2}}e^{-\sigma_0v_1D^2}||\\
&&\cdots ||A_{[+]}(1+D^2)^{-\frac{1}{2}}||||(1+D^2)^{\frac{1}{2}}e^{-\sigma_0v_mD^2}||dv\\
&&
\cdot \int_{\triangle_l}||e^{-\sigma_1v'_0D^2}||||A_{[+]}(1+D^2)^{-\frac{1}{2}}||||(1+D^2)^{\frac{1}{2}}e^{-\sigma_1v'_1D^2}||\\
&&\cdots ||A_{[+]}(1+D^2)^{-\frac{1}{2}}||||(1+D^2)^{\frac{1}{2}}e^{-\sigma_1v'_lD^2}||dv'\\
&&\leq \delta_0||C_0||||C_1||\sum_{l.m}\int_{\triangle_1}(\sigma_0)^{\frac{m}{2}}(\sigma_1)^{\frac{l}{2}}
\int_{\triangle_m}\int_{\triangle_l}\\
&&\cdot||\delta_1A_{[+]}(1+D^2)^{-\frac{1}{2}}||^{m+l}
{v}^{-\frac{1}{2}}_1\cdots {v}^{-\frac{1}{2}}_m
{v'}^{-\frac{1}{2}}_1\cdots {v'}^{-\frac{1}{2}}_ldvdv',\\
&&\leq \delta_0||C_0||||C_1||\sum_{l.m}\int_{\triangle_m}\int_{\triangle_l}\\
&&\cdot||\delta_1A_{[+]}(1+D^2)^{-\frac{1}{2}}||^{m+l}
{v}^{-\frac{1}{2}}_1\cdots {v}^{-\frac{1}{2}}_m
{v'}^{-\frac{1}{2}}_1\cdots {v'}^{-\frac{1}{2}}_ldvdv',\\
&&\leq\delta_0||C_0||||C_1||\left[\sum_{m, {\rm even}}||\delta_1A_{[+]}(1+D^2)^{-\frac{1}{2}}||^{m}\frac{2\pi^{\frac{m}{2}}}{(\frac{m+1}{2})!}\right.\\
&&\left.+\sum_{m, {\rm odd}}||\delta_1A_{[+]}(1+D^2)^{-\frac{1}{2}}||^{m}\frac{2\pi^{\frac{m}{2}}}{(\frac{m+1}{2})!}\right]^2\\
&&\leq\delta_0||C_0||||C_1||(1+\delta_2)^2e^{2(||\delta_1A_{[+]}(1+D^2)^{-\frac{1}{2}}||)^2}
~~~~~~~~~~~~~~~~~~~~~~~(4.49)
\end{eqnarray*}
where $\delta_0,\delta_1,\delta_2$ are constant.
For the general $q$, similarly we get
$$\left|\int_{\triangle_q}{\rm ^bStr^{end}}\left[C_0e^{-\sigma_0F}C_1e^{-\sigma_1F}\cdots C_q
e^{-\sigma_qF}\right]d\sigma\right|$$
$$\leq \delta_0\frac{1}{q!}(q+1+{\rm dim}X)\left(\prod_{j=0}^q||C_j||\right)\left(\delta_1e^{2||\delta_1A_{[+]}(1+D^2)^{-\frac{1}{2}}||^2}\right)^{q+1}.\eqno(4.50)$$
In order to estimate the type II integral, we decompose the type II integral as (4.36). Up to the last term, other terms have the same estimate
with corresponding terms in (4.36). Using the same trick as in [Xi], we get that the bound of the $1$-norm of the last term is
$\delta\frac{1}{q!}(\delta')^{q+1}(q+1)||B||\prod_{j=1}^q||A_j||$.\\
\indent By the above estimate,
${\rm ^bch}(B)$ is well-defined. Similarly, for a fixed $t>0$, ${\rm ^bch}(B_t)$ and ${\rm ^bch}(B_t,\frac{dB_t}{dt})$ are well defined. $\Box$\\

\section{The family Atiyah-Patodi-Singer index theorem for twisted Dirac operators}

\quad In this section, we extend the Getzler's index theorem to the family case. Let
$$\widehat{A}(R^{{\widehat{N}}/X})={\rm det}^{\frac{1}{2}}\left(\frac{\frac{R^{{\widehat{N}}/X}}{4\pi\sqrt{-1}}}{{\rm sinh}\frac{R^{{\widehat{N}}/X}}{4\pi\sqrt{-1}}}\right).\eqno(5.1)$$\\

\noindent {\bf Theorem 5.1} ~{\it  Suppose that all $D_{M,z}$ are invertible
with $\lambda$ the smallest positive eigenvalue of all $|D_{M,z}|$. We assume that
$||d(p|_M)||<\lambda$ and $p\in M_{r\times r}(C^{\infty}_{\rm exp}(\widehat{N}))$, then in the cohomology of $X$}
$${\rm ch}[{\rm Ind}(pD_{z,+}p)]=\int^b_{\widehat{N}/X}\widehat{A}(R^{{\widehat{N}}/X}){\rm ch}({\rm Imp})-\langle
{\widehat{\eta}}^*(B^M), {\rm ch}_*(p_M)\rangle.\eqno(5.2)$$

\noindent{\bf Proof.} By Theorem 4.6 and $(B+b)({\rm ch}(p))=0$, for fixed $t_1,t_2>0$, we have in the cohomology of $X$,
$$\langle{\rm ^bch}^*(B_{t_2}),{\rm ch}_*(p)\rangle-\langle{\rm ^bch}^*(B_{t_1}),{\rm ch}_*(p)\rangle=
-\frac{1}{\sqrt{\pi}}\langle\int^{t_2}_{t_1}{\rm ch}^*(B^M_t,\frac{dB^M_t}{dt})dt,{\rm ch}_*(p_M)\rangle.\eqno(5.3)$$
Let $t_1$ go to zero and $t_2$ go to $+\infty$. By Proposition 5.2 and Theorem 5.3 in [LMP], similar to the computations in Section 4 in [Wa3],
we get
$${\rm lim}_{t\rightarrow 0}{\rm ^bch}^{2k}(B_{t})(a_0,a_1,\cdots,a_{2k})=\frac{1}{(2k)!}(2\pi\sqrt{-1})^{-\frac{n}{2}}$$
$$\cdot\int^b_{\widehat{N}/X}
a_0da_1\wedge\cdots\wedge da_{2k}\widehat{A}(2\pi\sqrt{-1}R^{{\widehat{N}}/X}).\eqno(5.4)$$
Then we have
$${\rm lim}_{t_1\rightarrow 0}\langle{\rm ^bch}^*(B_{t_1}),{\rm ch}_*(p)\rangle=\int^b_{\widehat{N}/X}\widehat{A}(R^{{\widehat{N}}/X}){\rm ch}({\rm Imp}).\eqno(5.5)$$
By Lemma 5.2 in the following, we have
$${\rm lim}_{t_2\rightarrow +\infty}\langle{\rm ^bch}^*(B_{t_2}),{\rm ch}_*(p)\rangle={\rm lim}_{t\rightarrow +\infty}{\rm ^bch}^*(pB_{t}p).\eqno(5.6)$$
By all $D_{M,z}$ being invertible and Proposition 15 in [MP1], we have
$${\rm ch}[{\rm Ind}(pD_{z,+}p)]={\rm lim}_{t\rightarrow +\infty}{\rm ^bch}^*(pB_{t}p),\eqno(5.7)$$
By (5.3) and (5.5)-(5.7) and the definition of the eta cochain form, we get Theorem 5.1.~~$\Box$\\

\noindent {\bf Lemma 5.2} ~{\it The formula (5.6) holds.}\\

\noindent{\bf Proof.} Let $B_{t,u}=\sqrt{t}\psi_t(B+u(2p-1)[B,p])$. Using the same discussions with Theorem 4.6, we get in the cohomology of $X$
$$
\langle\frac{\partial{\rm ^bch}^*(B_{t,u})}{\partial u},{\rm ch}_*p\rangle
=-\frac{1}{\sqrt{\pi}}\langle{\rm ch}^*(B^M_{t,u},\frac{\partial B^M_{t,u}}{\partial u}),{\rm ch}_*(p_M)\rangle
.\eqno(5.8)$$
Then
$$\langle{\rm ^bch}^*(B_{t,1}),{\rm ch}_*p\rangle-\langle{\rm ^bch}^*(B_{t}),{\rm ch}_*p\rangle
=-\frac{1}{\sqrt{\pi}}\langle\int^1_0{\rm ch}^*(B^M_{t,u},\frac{\partial B^M_{t,u}}{\partial u})du,{\rm ch}_*(p_M)\rangle
.\eqno(5.9)$$
By $[B_{t,1},p]=0$, it holds that
$$\langle{\rm ^bch}^*(B_{t,1}),{\rm ch}_*p\rangle={\rm ^bch}(pB_{t}p).\eqno(5.10)$$
By (5.9) and (5.10) and the following lemma, we know that Lemma 5.2 is correct.~~$\Box$\\

\noindent{\bf Lemma 5.3}~{\it The following equality holds}
$${\rm lim}_{t\rightarrow +\infty}\langle\int^1_0{\rm ch}^*(B^M_{t,u},\frac{\partial B^M_{t,u}}{\partial u})du,{\rm ch}_*(p_M)\rangle=0.\eqno(5.11)$$

 \noindent{\bf Proof.} By $[B^M_{t,u},p_M]=(1-u)[B^M_t,p_M]$ and $\frac{\partial B_{t,u}}{\partial u}=(2p-1)[B_t,p]$, we have
$$ \langle\int^1_0{\rm ch}^*(B^M_{t,u},\frac{\partial B^M_{t,u}}{\partial u})du,{\rm ch}_*(p_M)\rangle
=\sum^{+\infty}_{l=0}\frac{(2l)!}{l!}t^{l+\frac{1}{2}}\sum_{j=0}^{2l}(-1)^{j+l}$$
$$\cdot\int_{u\in[0,1]}(1-u)^{2l}\psi_t\langle
p_M-\frac{1}{2},[B^M,p_M],\cdots,(2p_M-1)[B^M,p_M],\cdots,[B^M,p_M]\rangle_{t,u}.\eqno(5.12)$$
For the large $t$, we have
 $$||{\rm tr}e^{-tD^{M,2}_u}||\leq c_0e^{-t(\lambda-u|dp_M|)^2}.\eqno(5.13)$$
In the following, we drop off the index $M$. Using the same trick in Lemma 4.2 in [Wa3] and (5.13), we get the following estimate.
 For any $1\geq \sigma>0$, $t>0$ and $t$ is large
and any order $l$ fibrewise
differential operator ${A}$ with form coefficients, we have
$$||e^{-\sigma tB^2_u}A||_{\sigma^{-1}}\leq C_0(\sigma t)^{-\frac{l}{2}+\frac{{\rm dim}X}{2}}e^{-[(1-\varepsilon)(\lambda-u|dp|)^2-\varepsilon]\sigma t},\eqno(5.14)$$
where $C_0$ is a constant and $\varepsilon$ is any small positive constant.
By (5.14) and the H\"{o}lder inequality, we have
$$\left|\langle
p-\frac{1}{2},[B,p],\cdots,(2p-1)[B,p],\cdots,[B,p]\rangle_{t,u}\right|$$
$$\leq C_0\frac{|[B,p]||^{2l+1}}{(2l+1)!}t^{\frac{{\rm dim}X}{2}}e^{-[(1-\varepsilon)(\lambda-u|dp|)^2-\varepsilon]t}.\eqno(5.15)$$
By (5.12) and (5.15), we get
$$ \langle\int^1_0{\rm ch}^*(B^M_{t,u},\frac{\partial B^M_{t,u}}{\partial u})du,{\rm ch}_*(p_M)\rangle=O(e^{c_0(||dp||-\lambda)t}),\eqno(5.16)$$
where $c_0$ is a positive constant, so Lemma 5.3 holds.~~$\Box$\\

\noindent{\bf Acknowledgment.}~~This work
was supported by NSFC No.11271062 and NCET-13-0721. The author would like to thank Profs.
Weiping Zhang and Huitao Feng for introducing index theory to
him.  The author would
like to thank the referee for his careful reading and helpful comments.\\

\noindent{\bf Reference}\\

 \noindent [APS] M. F. Atiyah, V. K. Patodi and
I. M. Singer, {\it Spectral asymmetry and Riemannian geometry},
Math. Proc. Cambridge Philos. Soc. 77 (1975), 43-69; 78 (1975),
405-432; 79 (1976), 71-99.

\noindent[BeC]M. Benameur and A. Carey, {\it Higher spectral flow and an entire bivariant JLO cocycle}, J. K-theory, 11 (2013), 183-232.

\noindent [BGV] N. Berline, E. Getzler, and M. Vergne, {\it Heat
kernels and Dirac operators}, Spring-Verlag, Berline Heidelberg,
1992.

\noindent[BC] J. M. Bismut and J. Cheeger, {\it Families index theorem for manifolds with boundary
I, II,} J. Functional Analysis, 89 (1990), pp 313-363 and 90 (1990), pp 306-
354.

\noindent [BF] J. M. Bismut and D. S. Freed, {\it The analysis of
elliptic families II}, Commun. Math. Phys. 107 (1986), 103-163.

\noindent [CH] S. Chern and X. Hu, {\it Equivariant Chern character
for the invariant Dirac operators}, Michigan Math. J. 44 (1997),
451-473.

\noindent [Co]A. Connes, {\it Noncommutative differential geometry.} IHES. Publ. Math. No. 62 (1985), 257-360.

\noindent [CM]A. Connes and H. Moscovici, Transgression and
Chern character of finite dimensional K-cycles,
 Commun. Math. Phys. 155 (1993), 103-122.

\noindent [Do] H. Donnelly, {\it Eta invariants for G-space}, Indiana
Univ. Math. J. 27 (1978), 889-918.

\noindent [Fe] H. Feng, {\it A note on the noncommutative Chern
character (in Chinese)}, Acta Math. Sinica 46 (2003), 57-64.

\noindent [Ge1]E. Getzler, {\it Pseudodifferential operators on
supermanifolds and the Atiyah-Singer index theorem}, Commun. Math.
Phys. 92 (1983), 163-178.

\noindent [Ge2] E. Getzler, {\it Cyclic homology and the
Atiyah-Patodi-Singer index theorem}, Contemp. Math. 148 (1993),
19-45.

\noindent [GS] E. Getzler, and A. Szenes, {\it On the Chern
character of theta-summable Fredholm modules}, J. Func. Anal. 84
(1989), 343-357.

\noindent [LYZ] J. D. Lafferty, Y. L. Yu and W. P. Zhang, {\it A
direct geometric proof of Lefschetz fixed point formulas}, Trans.
AMS. 329 (1992), 571-583.

\noindent[LMP]M. Lesch, H. Moscovici and J. Pflaum, {\it Connes-Chern character for manifolds with boundary and eta cochains}, Mem. AMS 220 (2012), no.1036, viii+92 pp.

 \noindent[MP1]R. Melrose and P. Piazza, {\it Families of Dirac operators, boundaries and the b-calculus}, J. Differential Geom. 46 (1997), no. 1, 99-180.

 \noindent[MP2]R. Melrose and P. Piazza, {\it An index theorem for families of Dirac operators on odd-dimensional manifolds with boundary}, J. Differential Geom. 46 (1997), no. 2, 287-334.

 \noindent[Po]R. Ponge, {\it A new short proof of the local index
formula and some of its applications}, Comm. Math. Phys. 241(2003), 215-234.

\noindent[PW]R. Ponge and H. Wang, {\it Noncommutative geometry, conformal geometry, and the local equivariant index theorem},
 J. Noncommut. Geom. 10 (2016), 307-378.

\noindent[Wa1]Y. Wang, {\it The equivariant noncommutative Atiyah-Patodi-Singer index theorem}, K-Theory 37 (2006),
213-233.

\noindent[Wa2]Y. Wang, {\it The noncommutative infinitesimal equivariant index formula}, J. K-Theory 14 (2014), 73-102,

\noindent[Wa3]Y. Wang, {\it Volterra calculus, local equivariant family index theorem and equivariant eta forms}, arXiv:1304.7354 to appear in Asian J. Math.

\noindent[Wa4]Y. Wang, {\it The noncommutative infinitesimal equivariant index formula: part II}, J. Noncommut. Geom.  10 (2016), 379-404.

\noindent[Xi]Z. Xie, {\it The odd-dimensional analogue of a theorem of Getzler and Wu}, J. Noncommut. Geom. 7 (2013), no. 3, 647-676.

\noindent [Wu] F. Wu, {\it The Chern-Connes character for the Dirac
operators on manifolds with boundary}, K-Theory 7 (1993), 145-174.

\noindent [Zh] W. P. Zhang, {\it A note on equivariant eta
invariants}, Proc. AMS. 108 (1990), 1121-1129.\\

 \indent{\it School of Mathematics and Statistics, Northeast Normal University, Changchun Jilin, 130024, China }\\
 \indent E-mail: {\it wangy581@nenu.edu.cn}\\

\end{document}